\documentclass[a4paper]{siamltex}

\usepackage{a4wide}
\usepackage[english]{babel}
\usepackage{amssymb,amsfonts,amsmath,latexsym}
\usepackage{graphicx}
\usepackage{url}

\newcommand{\R}{{\mathbb R}}

\numberwithin{equation}{section}

\begin{document}

\title{Behavior patterns of 15-year-old students in a Digital
Storytelling in Mathematics. A social network analysis}

\author{A.Concas\thanks{Department of Mathematics and Computer Science,
University of di Cagliari, via Ospedale 72, 09123 Cagliari, Italy. E-mail:
\texttt{anna.concas@unica.it}, \texttt{mpolo@unica.it}},
\and
M. Polo\footnotemark[1]}

\maketitle


\begin{abstract}

This paper presents the preliminary considerations of the application of a software to an experimental work conducted on the Digital
Storytelling in Mathematics, as part of the project Prin 2015 ``Digital Interactive Storytelling in Mathematics: a competence-based social approach''. 
An activity designed for promoting the critical mathematical thinking among the students that foresees them to participate as active protagonists and as observers of the protagonists during the problem solving activity, will be illustrated and then the outcomes will be examined from a numerical analysis point of view. In particular, the interactions between the partecipants will be investigated by using a Matlab software for solving the \emph{seriation problem}.
\end{abstract}


\section{Introduction}
In this paper, the e-learning topic in the mathematical education is addressed from a point of view of an analysis of the potentiality and also of the drawbacks of using online platforms as an instrument for helping the learning process.
The aim consists in proving how, also in this context, virtual interactions between the participants is fundamental.
Indeed it can be observed that the interactions within groups of students or with an expert may modify the educational development.

In particular, we present some results of an analysis of the social interactions between the students of a class which took part to a trial project in January 2020 within the PRIN project titled ``Digital Interactive Storytelling in
Mathematics: A Competence-based Social Approach''.

This work refers to an experimentation related to a design online competence-
oriented activities, making use of narrative and social approaches, developed within the
aforementioned PRIN research project.

In this paper we will focus on the analysis of a case study  concerning the analysis, by using numerical techniques in the the complex network theory framework, of the behaviours of the teacher and the students in an activity devoted to promote the critical thinking and argumentation in mathematics \cite{APP19, APP20}. The aim of the work is to validate the toolbox~\cite{CFR} in a didactis framework and shed light on how the analysis
of complex networks might be fruitful in a DIST-M (Digital Interactive
Storytelling in Mathematics) environment.

\section{Theoretical framework}
The whole PRIN project is framed in a Vygotskian approach in which social interaction and communication play a key role in the evolution of students and language has been widely recognized as fundamental in learning mathematics. Students' participation in social interactions and in communication is the key to their evolution and the appropriation of cultural tools, such as language; see~\cite{vyg12}. According to Sfard~\cite{sfard}, thought is a form of communication and languages do not are just couriers of pre-existing meanings, but they are builders of meanings themselves.

The idea of implementing a teaching/learning process into a digital interactive storytelling  using Moodle, where students are immersed as actors or observers, fits the idea of expanding the context of mathematics education beyond the classroom~\cite{En17}. This expansion occurs along the layer of bridging elements from societal practices in mathematics instruction: indeed, students are more and more involved in digital games/stories outside the school. The DIST-M would bring such elements into mathematics education practices, bridging the agency of experts in mathematics knowledge and the motivation of the students. It also allows to add on the school activity system with new societal components, such as rules, community and division of labour~\cite{laz16}.
In the presented case study, students and teacher are engaged in a suitable designed online activity, where the interactions among students and students-teacher are well organized. The organization of the activity foresees that the roles can be played as actors or as observers ~\cite{APP20}. Nevertheless, ``the very subject of learning is transformed from isolated individuals to collectives''~\cite{ENSAN}. Thus learning does not depend on a single role and the related actions, but on the whole group, as a system. In this paper we will analyze the outcome of a pilot, with particular reference to the students who played the roles played as observers.

Following Wegner~\cite{weg99} the learning mechanisms as negotiation, 
as a response to the pedagogical intentions of the context, but also considering that the purpose of educational planning is not to appropriate learning and institutionalize it in a structured process, but to support the formation of learning communities. 

The learning process is interpreted as a social way of changing the participation from peripheral to a community practice involvement.
This process develops in a cooperative approach ~\cite{pesci09, sulisworo2016cooperative} which aims to create a situation where individual
success is determined or influenced by the success of the group.

This learning view point and the collaborative learning~\cite{pesci09} complement each other even if the intercation roles of equals are not a priori established.

The roles assigned to the actors who participated to the project, allowed 	the cohesiveness between the knowledges leading to a real educational organization~\cite{chev07} and the cohesion between the the different perspectives of the participants plays a crucial function~\cite{chev08}.

In Polo et all.~\cite{PIFP,APP19} the authors presented a social analysis of the interactions among the students involved in a trial of the PRIN project ``Digital Interactive Storytelling in Mathematics: a Competence-based Social Approach''. The project is designed within a digital storytelling framework where the story follows the interactions among the characters played by the students and an expert (teacher or researcher). We report the results of a trial that involved teachers and students from the upper secondary school analysing from a Social Network Analysis point of view the interventions of the expert, the involvement / participation of the students and the interactions among peers and with the expert. We also briefly discuss the potential and limitations of the currently available tools to perform this kind of analysis, in view of the much broader perspective offered by the Learning Analytics approach.

As mentioned previously, this work regards a social network analysis carried out by modelling the experimentation trial as a seriation problem. Online communication tools contribute the creation of communities of practice and learning online allowing the interactions in an online education environment; see ~\cite{FR18} and the references therein. In this framework, the seriation problem consists in finding an ordering of the set of considered units, which represent the students in a specific group according to a given correlation function so that elements with higher similarity in the behaviour  are close to each other in the resulting sequence obtained by applying a toolbox recently developed in Matlab~\cite{CFR}
\medskip

%
%

\section{Methodology and research hypothesis}
The activities are framed in a narrative, in which there are characters that will act as avatars for the various participating students.  For more details about the design of the activity see~\cite{APP19,APP20}. The activites included in the trial can be contextualize in the Digital Interactive Storytelling framework described in~\cite{PIFP}.
\subsection{Didactical modeling}

The design is structured in such a way that it includes both individual and collaborative work and discussion with the expert who acts as a mediator by intervening when necessary. This moderation concerns not only helping the students with the mathematical problem they have to solve but it affects also the communicative part when the students have to formalize the solution they found together. Different tools have been considered in the analysis.  Communication tools, like the chat the students used to communicate within their group, through which we were able to observe the usage of different kind of registers (colloquial vs. literate). Tools like shared files allowed the students to develop a higher collaboration; while supporting tools such as spreadsheets, CAS, calculators and word blocks have been made available for inquiry, formalization and proof.


The schema for the entire DIST-M activity evolves into different phases, each one corresponding to a specific episode in the story.
The first phase is the exploration one where a concise description of what has been
observed is proceduced. The description is then refined in the second phase in order to obtain a statement that
serves as a conjecture. This utterance, is manipulated in the third phase in order to arrive at its appropriate
formalization  and to be adjusted to lead to a
proof. In the fourth stage, students are asked to identify and organize the arguments
into an appropriate deductive chain, justifying each step of the deduction, in order to construct the proof.
In the last phase, reflection is required, both as an evaluation on the cognitive level -students have to
read again the story from the mathematical point of view- and as a self-assessment, on the metacognitiveand affective level -each student is asked to think about the roles played as an ``actor'' or as an ``active
observer''.

The instructional activities take place in the context of a narrative, in a situation that aims to be engaging
and familiar to the student. The story evolves over time and each student is a character in the story
(Albano et al., 2019). Personal and group interaction is moderated by the teacher playing herself a
character in the story. As the narrative progresses, the story evolves in accordance with the characters'
interaction with it.
The genre chosen was that of science fiction, which sees a group of four friends Marco, Clara, Federico
and Sofia engaged in the challenge of communicating with the aliens from whom they had received
some mysterious messages made of numbers and operations. Alongside the friends, there is also an
adult, who is the avatar of the teacher/tutor and acts as an expert in the learning path.
The story is organized in 5 episodes, where students play as characters: Marco (the Boss) is the leader
of the group and takes care to organize the tasks among his peers; Clara (the Pest) intervenes by posing
doubts and question about the mathematical problem; Sofia (the Blogger) loves writing and takes care
to summarizes the answers shared by the group; Federico (the Promoter) proposes new ideas and paths
to follow. They are given a role (which changes from one episode to the other) and some actions to
perform, according to the role the students are playing in that episode. Some tasks are individual, some
others are collaborative and therefore require a good grade of communication and coordination
(students can talk to each other in a chat). In order to allow all students to play within the story, we
arranged them in teams of 4 members (one for each of the above characters). In each episode, only one
team is actively involved within the story and each member plays a different character; they can interact
among themselves and with the story. In the meantime, the members of the other teams play as
Observers, each member of a team observing a different character. These are actively involved in taking
notes in their personal logbooks, reporting how (and how well) the observed active character performed.
In the subsequent episode, another team becomes active and all the roles are switched in turn. In this
way all students will experience all roles and all activities of the story, though not always as an active
player.


\subsection{Numerical modeling}

From a numerical point of view, this paper is related to an application of seriation to the Social Network Analysis and in particular for analysing the roles played by students in the digital interactive storytelling pilot project described above.
Seriation is a fundamental combinatorial ordering problem, asking to find
the best enumeration order of a set of units, according to a given correlation function, so that
elements with higher similarity are close to each other in the resulting sequence. The desired
order can be characteristic of the data, a chronological order, a gradient or any sequential
structure of the data. We will state the seriation problem
from the mathematical point of view by considering it as the arrangement of units in a
sequence according to a ``similarity'' function.
The concept of seriation has been formulated in many different ways and the seriation
problem, along with seriation methods and algorithms, emerges and finds application in a wide range of contexts spanning from archaeology,
anthropology, to genomics and DNA sequencing, and, within mathematics, from complexity
theory to combinatorial optimisation, graph theory and operational research. The interested reader can find an overview on the context and application of seriation in~\cite{C20} and references therein.
The concept of seriation first arose in the archaeological field; see ~\cite{PAS} for a review on the seriation problem in
archaeology with a mathematical perspective.
In all the applications, seriation data are usually given in terms of a matrix, called \emph{data
matrix}, whose rows or columns (or both) represent the elements to be ordered. Considering
the application in this research project, the rows of the data matrix represent the partecipants while the columns represent the features divided, as explained below in details, into principal and secondary features. Each player is characterized by the role she/he plays in the considered episodes and by the partecipation in terms of messages written on chat.

The data matrix will be hereafter named \emph{adjacency matrix} since it describes a bipartite graphs  which allows to define the interrelationship between the units/roles we want to rearrange according to the observed variables.

Bipartite graphs, also called ``two-mode'' networks in the sociology literature, represent the membership of nodes in groups. Indeed, a graph is bipartite if its vertices can be split into two disjoint subsets such that only edges between nodes belonging to different sets can occur.
Formally, bipartite networks can be defined as follows.
\begin{definition}
A graph $G$ is \emph{bipartite} if its vertices can be divided into two disjoint sets $U$ and $V$ containing $n$ and $m$ nodes, respectively, such that
every edge connects a node in $U$ to one in $V$.
\end{definition}

In our research environment metaphor, the two disjoint nodes sets $U$ and $V$ represent the players (students and teacher) (i.e. units we want to somehow rearrange according to the role they play) and the features related to them, respectively.

In this framework, the adjacency matrix of order $n\times m$ can be interpreted as the upper-right block of the adjacency matrix of the associated  bipartite graph
$$
\mathcal{A_B} = 
\begin{bmatrix}
0_n & A \\ 
A^T & 0_m
\end{bmatrix}\in \R^{(n+m)\times (n+m)},
$$
with $0_n$  and $0_m$ null matrices of order $n$ and $m$ respectively and 
the adjacency matrix $A\in\R^{n\times m}$ describes the connections between the two node sets and it is obtained by setting 
\[a_{i,j}=
\begin{cases}
1 \quad \text{\footnotesize{if player $i$ presents features of type $j$}} \\
0 \quad \text{\footnotesize{otherwise}}
\end{cases}.
\]
In case of elements  $a_{i,j}\neq 1$, then the adjacency matrix is usually called abundance matrix in the seriation problem. In this case every entry represents the relative frequency or the number of features (played role and considered variables) related to the considered player.
As stated above, we will follow the typical terminology used in complex networks theory and we refer to a binary representation of seriation data as an adjacency matrix.

The first mathematical
definition of seriation was based on the construction of a symmetric matrix $S$ known as
\emph{similarity} matrix,  where the element $s_{i,j}$ describes, in some way, the likeness of
the nodes $i, j \in U$ representing two units (i.e. partecipants). One possible definition of the
similarity matrix is through the product $S = AA^T$, being $A$ the adjacency matrix of the considered 
problem. Even in the presence of abundance data, the matrix we will consider as adjacency matrix describing the case study will be binary and constructed by setting equal to 1 every non zero entry. This choice is required by the fact that the software used for processing the data can only, up until now, manipulate binary seriation data.
In the similarity matrix the largest value on each row is the diagonal
element, which reports the number of variables associated to each unit/partecipant. By permuting the rows
and columns of $S$ in order to cluster the largest values close to the main diagonal, one obtains
a permutation of the corresponding rows of $A$ that places the people similar in behaviour closer to
each other. 
Summarizing, one can map the set of pairwise relative measurements among the objects
to a similarity matrix which represents the information of the objects/units to be ordered. Then,
the seriation problem can be modeled as a discrete optimization problem, whose goal is
to find a permutation of the rows and columns of the similarity matrix which minimizes a
given objective function. Methods based on the similarity matrix differ from one to another
on how the similarity matrix is constructed or because of different objective functions used
for evaluate the found permutations of the seriation data matrix.

In this paper we applied a Matlab toolbox described in~\cite{CFR}  with the aim of validating the software in a didactical research setting. The toolbox contains an implementation of a spectral algorithm for obtaining the solution
of the seriation problem based on the use of an
eigenvector, called Fiedler vector and associated to the eigenvalue, called Fiedler value of
the Laplacian matrix associated to the problem. 
The algorithm, by using the sorted entries of the Fiedler
vector, returns a classification of the ordering permutations in terms of a compact data structure called
PQ-tree briefly described hereafter. The package, named PQser toolbox, also defines a data structure to store a
PQ-tree and provides the Matlab functions to manipulate and visualize it. The detailed description of
the spectral algorithm is skipped here for the sake of brevity, but the interested reader can refer to~\cite{CFR}. 
In the following subsection we focuse on the outcome of the spectral algorithm that encode all the possible reorderings of the units/roles.

A PQ-tree is a data structure introduced by Booth and Lueker~\cite{BL76} to encode a family
of permutations of a set of elements and solve problems connected to finding admissible
permutations according to specific rules.
A PQ-tree $T$ over a set $U = \{u_1 , u_2 , \ldots, u_n\}$ is a rooted, ordered tree whose leaves are
in one-to-one correspondence with the elements of $U$ and their order gives a reordering of
the elements of the considered set. The internal (non-leaf) nodes of $T$ are distinguished
as either P-nodes or Q-nodes. The only difference between them is the way in which
their children can be reordered. Namely, for a P-node, its children may be
arbitrarily reordered so all possible permutations of the children leaves are permitted; for a
Q-node only one order and its reverse are allowed since the children leaves may be ordered
only left-to-right or right-to-left. The root of the tree can be either a P or a Q-node.
We will represent graphically a P-node by a circle, and a Q-node by a rectangle. The leaves of T will be
displayed as triangles, and labeled by the elements of U. The frontier of T is one possible
permutation of the elements of U, obtained by reading the labels of the leaves from left to
right.

\begin{figure}
\begin{minipage}{.52\textwidth}
\includegraphics[width=\textwidth]{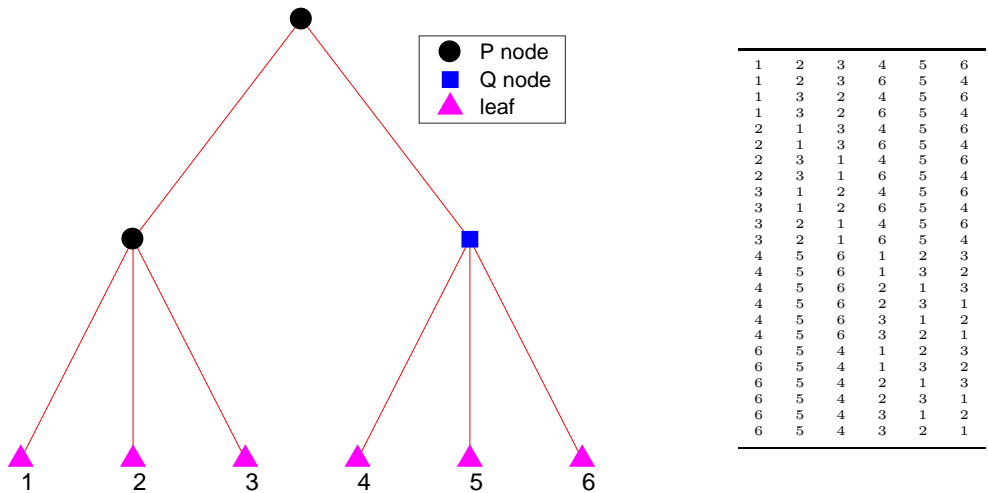}
\end{minipage}
\begin{minipage}{.45\textwidth}
\tiny
\begin{center}
\begin{tabular}{cccccc}
\hline\noalign{\smallskip}
1 & 2 & 3 & 4 & 5 & 6 \\
1 & 2 & 3 & 6 & 5 & 4 \\
1 & 3 & 2 & 4 & 5 & 6 \\
1 & 3 & 2 & 6 & 5 & 4 \\
2 & 1 & 3 & 4 & 5 & 6 \\
2 & 1 & 3 & 6 & 5 & 4 \\
2 & 3 & 1 & 4 & 5 & 6 \\
2 & 3 & 1 & 6 & 5 & 4 \\
3 & 1 & 2 & 4 & 5 & 6 \\
3 & 1 & 2 & 6 & 5 & 4 \\
3 & 2 & 1 & 4 & 5 & 6 \\
3 & 2 & 1 & 6 & 5 & 4 \\
4 & 5 & 6 & 1 & 2 & 3 \\
4 & 5 & 6 & 1 & 3 & 2 \\
4 & 5 & 6 & 2 & 1 & 3 \\
4 & 5 & 6 & 2 & 3 & 1 \\
4 & 5 & 6 & 3 & 1 & 2 \\
4 & 5 & 6 & 3 & 2 & 1 \\
6 & 5 & 4 & 1 & 2 & 3 \\
6 & 5 & 4 & 1 & 3 & 2 \\
6 & 5 & 4 & 2 & 1 & 3 \\
6 & 5 & 4 & 2 & 3 & 1 \\
6 & 5 & 4 & 3 & 1 & 2 \\
6 & 5 & 4 & 3 & 2 & 1 \\
\noalign{\smallskip}\hline
\end{tabular}
\end{center}
\end{minipage}\smallskip
\caption[PQ-tree over the set $U=\{1,\ldots,6\}$ and the encoded permutations.]{On the left, a PQ-tree over the set $U=\{1,\ldots,6\}$; on the right,
the 24 admissible permutations encoded in the tree.}
\label{fig:pqtree}
\end{figure}

Hence, a PQ-tree represents in a compact way permutations of the elements of a set
through admissible reorderings of its leaves. Each transformation specifies
an allowed reordering of the nodes within a PQ-tree. For example, a tree with a single
P-node represents the equivalence class of all permutations of the elements of U, while a
tree with a single Q-node represents both the left-to-right and right-to-left orderings of the
leaves. A tree with a mixed P-node and Q-node structure represents the equivalence class of
a constrained permutation, where the structure of the tree determines the constraints on the
admissible permutations.
Figure~\ref{fig:pqtree} displays a PQ-tree with a P node as root and all the admissible permutations it represents and encodes.
Both the figure and the permutations encoded in the tree have been obtained by applying
specific functions contained in the PQser toolbox.

Now, we briefly decribe the spectral algorithm we applied for analizing the case study as described in the following. 
In case of a reducible similarity matrix, it can be proved that the seriation problem can be decoupled~\cite{atk}. The spectral algorithm detect irreducible blocks of the similarity matrix identifying also the
corresponding index sets. The irreducible blocks correspond to the connected components of the considered graph. If more than one connected component is found, then the function calls itself on
each component, and stores the returned lists of nodes as children of a P-node. If the similarity matrix of the considered problem is irreducible, the dimension $n$ of the matrix is considered. The cases $n = 1,\quad 2$ are trivial. If $n > 2$, the Laplacian matrix $L$ of $S$ is computed, as well as the Fiedler value and the correspinding eigenvector. The PQser toolbox computes by default the eigenpairs corresponding to the
three eigenvalues of smallest magnitude, since they are sufficient to understand if the
Fiedler value is simple or multiple, but the default value can be modified. The algorithm determines the multiplicity of the Fiedler value according to a given tolerance. If the Fiedler value is a simple eigenvalue of L, the algorithm sorts the elements of the current list according to the reordering of the Fiedler vector and stores them as the children of a Q-node. If two or more values of the Fiedler vector
are repeated, the function invokes itself recursively; on the contrary, the corresponding node becomes a leaf In the presence of a multiple Fiedler vector, the problem is not well posed and an approximate solution is computed, as described in~\cite{CFR}.
The number $N$ of all the admissible permutations generated by the algorithm can
be obtained by counting all the admissible boundaries of the tree.

In this framework, with the aim of validating the PQser toolbox in a didatics field, we applied the software in order to find all the admissible reordering of the units, i.e., the people involved in the considered pilot project on the basis of their behaviour studied through the chat they used during the trial.



\section{Case study}
In this section we analyze the first episode of the final prototype experimented in January 2020.
The group of the participants of this project is composed by 26 students and the teacher who, as already mentioned, assumed the role of a mediator in the group of the 4 students who play the actors role.
The story evolves in 5 episodes and in each one of them the students play a different role being an actor or an observer in a different group; see ~\cite{APP19}.
The aim behind the narration is the solution of a mathematical problem by proving that, given four consecutive natural numbers, the difference between the product of the first and the fourth and the product between the second and the third is always equal to 2.
The four principal roles, assigned to the students who are actors, are of the \emph{Peste Clara}, the \emph{Blogger Sofia}, the \emph{Boss Marco } and the \emph{Promoter Federico}.

%

All the 26 students have been divided into 6 groups indicated thereafter as  $g_i$ for $i = 1,\ldots,6$. The group $g_1$ identifies, in each episode of the narration, the group of the actors in the considered framework of the story. The other groups, in which every member assumes one of the four above-mentioned principal roles, represent the observers. The task of the latter is that of observing the actors who interact with each other focusing, in particular, on the behaviour of the actor who assumes their same role.

The investigation has been accomplished by analysing the group chat into which the studends who act interacted with each other and with the teacher who helps the actors in their assignement. Also the chat used by the observers have been examined by using the Matlab sotware published in~\cite{CFR} and described above.
Related to the complex network theory, the chat examination led to the construction of a matrix whose rows represent all the 27 participants to the project while the columns report the principal and secondary variables consisting of the roles and the interactions divided into primary and minor too.

The adjacency matrix considered in the esperimentation, of size $27\times 31$, represent abundance data and it is reported below. In particular, the matrix is such that the rows and columns are costructed as follows.
The rows are 27 and contain all the partecipants to the trial: the first four represent the students who act in the considered episode. The rows from the 5th to the second-last are associated to the other students who assume the role of observers, while the last one is related to the teacher who mediates the actors participation in the chat.
The columns are instead 31 and they contain the roles on each group into which the students are divided and the variables considered in the pilot project. More in detail, the columns from the 1st to the 4th represent the roles assumed by the actors, the fifth represents the role taken on by the teacher. The columns from the 6th to the 25th represent the roles assumed by the students on each observers group.

The last columns represent the variables which take into account the interactions between the student and with the teacher. Precisely, column 26 report the interactions between the students related to the mathematical problem, column 27 is related to the social interactions between the actors. The 28th column reports the interaction of mathematical nature between the actors and the teacher; column 29 takes into account the social interactions between the students and the teacher. The 30th column reports the mathematical exchanges between students and teacher and in the last one the element $a_{i,31}$ represents the social interactions between the teacher and the $i$-th student.

$$\tiny{A=\left(\begin{array}{ccccccccccccccccccccccccccccccc} 1 & 0 & 0 & 0 & 0 & 0 & 0 & 0 & 0 & 0 & 0 & 0 & 0 & 0 & 0 & 0 & 0 & 0 & 0 & 0 & 0 & 0 & 0 & 0 & 0 & 18 & 8 & 7 & 6 & 5 & 2\\ 0 & 1 & 0 & 0 & 0 & 0 & 0 & 0 & 0 & 0 & 0 & 0 & 0 & 0 & 0 & 0 & 0 & 0 & 0 & 0 & 0 & 0 & 0 & 0 & 0 & 16 & 10 & 3 & 2 & 1 & 2\\ 0 & 0 & 1 & 0 & 0 & 0 & 0 & 0 & 0 & 0 & 0 & 0 & 0 & 0 & 0 & 0 & 0 & 0 & 0 & 0 & 0 & 0 & 0 & 0 & 0 & 16 & 6 & 5 & 2 & 3 & 0\\ 0 & 0 & 0 & 1 & 0 & 0 & 0 & 0 & 0 & 0 & 0 & 0 & 0 & 0 & 0 & 0 & 0 & 0 & 0 & 0 & 0 & 0 & 0 & 0 & 0 & 11 & 7 & 1 & 0 & 1 & 0\\ 0 & 0 & 0 & 0 & 0 & 1 & 0 & 0 & 0 & 0 & 0 & 0 & 0 & 0 & 0 & 0 & 0 & 0 & 0 & 0 & 0 & 0 & 0 & 0 & 0 & 0 & 0 & 0 & 0 & 0 & 0\\ 0 & 0 & 0 & 0 & 0 & 0 & 0 & 0 & 0 & 0 & 0 & 0 & 0 & 0 & 0 & 0 & 0 & 0 & 0 & 0 & 0 & 0 & 0 & 0 & 0 & 0 & 0 & 0 & 0 & 0 & 0\\ 0 & 0 & 0 & 0 & 0 & 0 & 0 & 1 & 0 & 0 & 0 & 0 & 0 & 0 & 0 & 0 & 0 & 0 & 0 & 0 & 0 & 0 & 0 & 0 & 0 & 0 & 0 & 0 & 0 & 0 & 0\\ 0 & 0 & 0 & 0 & 0 & 0 & 0 & 0 & 1 & 0 & 0 & 0 & 0 & 0 & 0 & 0 & 0 & 0 & 0 & 0 & 0 & 0 & 0 & 0 & 0 & 0 & 0 & 0 & 0 & 0 & 0\\ 0 & 0 & 0 & 0 & 0 & 0 & 0 & 0 & 0 & 1 & 0 & 0 & 0 & 0 & 0 & 0 & 0 & 0 & 0 & 0 & 0 & 0 & 0 & 0 & 0 & 0 & 0 & 0 & 0 & 0 & 0\\ 0 & 0 & 0 & 0 & 0 & 0 & 0 & 0 & 0 & 0 & 1 & 0 & 0 & 0 & 0 & 0 & 0 & 0 & 0 & 0 & 0 & 0 & 0 & 0 & 0 & 0 & 0 & 0 & 0 & 0 & 0\\ 0 & 0 & 0 & 0 & 0 & 0 & 0 & 0 & 0 & 0 & 0 & 1 & 0 & 0 & 0 & 0 & 0 & 0 & 0 & 0 & 0 & 0 & 0 & 0 & 0 & 0 & 0 & 0 & 0 & 0 & 0\\ 0 & 0 & 0 & 0 & 0 & 0 & 0 & 0 & 0 & 0 & 0 & 0 & 1 & 0 & 0 & 0 & 0 & 0 & 0 & 0 & 0 & 0 & 0 & 0 & 0 & 0 & 0 & 0 & 0 & 0 & 0\\ 0 & 0 & 0 & 0 & 0 & 0 & 0 & 0 & 0 & 0 & 0 & 0 & 0 & 1 & 0 & 0 & 0 & 0 & 0 & 0 & 0 & 0 & 0 & 0 & 0 & 0 & 0 & 0 & 0 & 0 & 0\\ 0 & 0 & 0 & 0 & 0 & 0 & 0 & 0 & 0 & 0 & 0 & 0 & 0 & 0 & 1 & 0 & 0 & 0 & 0 & 0 & 0 & 0 & 0 & 0 & 0 & 0 & 0 & 0 & 0 & 0 & 0\\ 0 & 0 & 0 & 0 & 0 & 0 & 0 & 0 & 0 & 0 & 0 & 0 & 0 & 0 & 0 & 1 & 0 & 0 & 0 & 0 & 0 & 0 & 0 & 0 & 0 & 0 & 0 & 0 & 0 & 0 & 0\\ 0 & 0 & 0 & 0 & 0 & 0 & 0 & 0 & 0 & 0 & 0 & 0 & 0 & 0 & 0 & 0 & 1 & 0 & 0 & 0 & 0 & 0 & 0 & 0 & 0 & 0 & 0 & 0 & 0 & 0 & 0\\ 0 & 0 & 0 & 0 & 0 & 0 & 0 & 0 & 0 & 0 & 0 & 0 & 0 & 0 & 0 & 0 & 0 & 1 & 0 & 0 & 0 & 0 & 0 & 0 & 0 & 0 & 0 & 0 & 0 & 0 & 0\\ 0 & 0 & 0 & 0 & 0 & 0 & 0 & 0 & 0 & 0 & 0 & 0 & 0 & 0 & 0 & 0 & 0 & 0 & 1 & 0 & 0 & 0 & 0 & 0 & 0 & 0 & 0 & 0 & 0 & 0 & 0\\ 0 & 0 & 0 & 0 & 0 & 0 & 0 & 0 & 0 & 0 & 0 & 0 & 0 & 0 & 0 & 0 & 0 & 0 & 0 & 1 & 0 & 0 & 0 & 0 & 0 & 0 & 0 & 0 & 0 & 0 & 0\\ 0 & 0 & 0 & 0 & 0 & 0 & 0 & 0 & 0 & 0 & 0 & 0 & 0 & 0 & 0 & 0 & 0 & 0 & 0 & 0 & 1 & 0 & 0 & 0 & 0 & 0 & 0 & 0 & 0 & 0 & 0\\ 0 & 0 & 0 & 0 & 0 & 0 & 0 & 0 & 0 & 0 & 0 & 0 & 0 & 0 & 0 & 0 & 0 & 0 & 0 & 0 & 0 & 0 & 0 & 0 & 0 & 0 & 0 & 0 & 0 & 0 & 0\\ 0 & 0 & 0 & 0 & 0 & 0 & 0 & 0 & 0 & 0 & 0 & 0 & 0 & 0 & 0 & 0 & 0 & 0 & 0 & 0 & 0 & 1 & 0 & 0 & 0 & 0 & 0 & 0 & 0 & 0 & 0\\ 0 & 0 & 0 & 0 & 0 & 0 & 0 & 0 & 0 & 0 & 0 & 0 & 0 & 0 & 0 & 0 & 0 & 0 & 0 & 0 & 0 & 1 & 0 & 0 & 0 & 0 & 0 & 0 & 0 & 0 & 0\\ 0 & 0 & 0 & 0 & 0 & 0 & 0 & 0 & 0 & 0 & 0 & 0 & 0 & 0 & 0 & 0 & 0 & 0 & 0 & 0 & 0 & 0 & 1 & 0 & 0 & 0 & 0 & 0 & 0 & 0 & 0\\ 0 & 0 & 0 & 0 & 0 & 0 & 0 & 0 & 0 & 0 & 0 & 0 & 0 & 0 & 0 & 0 & 0 & 0 & 0 & 0 & 0 & 0 & 0 & 1 & 0 & 0 & 0 & 0 & 0 & 0 & 0\\ 0 & 0 & 0 & 0 & 0 & 0 & 0 & 0 & 0 & 0 & 0 & 0 & 0 & 0 & 0 & 0 & 0 & 0 & 0 & 0 & 0 & 0 & 0 & 0 & 1 & 0 & 0 & 0 & 0 & 0 & 0\\ 0 & 0 & 0 & 0 & 1 & 0 & 0 & 0 & 0 & 0 & 0 & 0 & 0 & 0 & 0 & 0 & 0 & 0 & 0 & 0 & 0 & 0 & 0 & 0 & 0 & 0 & 0 & 0 & 0 & 14 & 6 \end{array}\right)}
$$

Let us consider the first student, by examining the first row of $A$, it can be observed that the student plays the role of the ``PESTE'' corresponding to the first column (the element $a_ {1,1} = 1 $ indicates that student 1 has assumed the first role which is that of the PESTE Clara). The student interacts with classmates 18 times (main interactions / on the mathematical problem) being the element $ a_ {1,26} = 18 $; element $ a_ {1,27} = 8 $ indicates that actor 1 carries out 8 interventions of social interaction with companions. The element $ a_ {1,28} = 7$ indicates that the PESTE student interacted 7 times with the teacher on the mathematical problem and the social / secondary interactions with the mediatior were instead 6 (element $ a_ {1,29} = $ 6). The element $ a_ {1,30} = 5 $ indicates that the teacher has interacted 5 times directly with the student on the mathematical problem, while the last column shows that the teahcre's social interactions with the student were 2.

By considering the data matrix $A$, it can be observed that the group of the actors given by the first 4 students interact only with the teacher. The other students, who are the observers in this phase of the experimentation, do not interact neither with each other nor with the first group and the teacher. This fact leads to the observation that, in terms of graph theory, the indices corresponding to the actors and the teacher may give a connected component of the underlying graph. 
Hence, by applying the spectral algorithm, the indices identitying the actors together and he teacher should correspond to the childrend of a Q node.

The other students, on the other hand, may be considered as isolated nodes and so their indices, according to how the spectral seriation algorithm works, should identify the children leaves of a P node.

\begin{figure}
\begin{center}
\includegraphics[width=.8\textwidth]{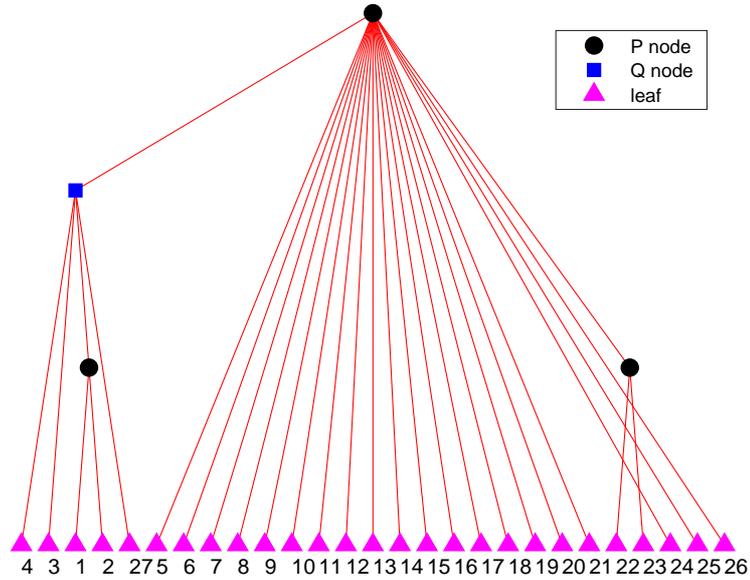}
\caption{PQ-tree corresponding to the data matrix $A$.}
\label{pqtree}
\end{center}
\end{figure}

By using the data matrix, A, we constructed the similarity matrix associated with our application.
The observation that the group of the actors given by the first 4 students interact only with the
teacher while the other students, who are the observers in this phase of the experimentation, interact
neither with each other nor with the first group and the teacher, leads to notice that, in terms of
graph theory, the indices corresponding to the actors and the teacher may give a connected
component of the underlying graph., i.e., set of vertices that are linked to each other and are
connected to no additional vertices in the rest of the graph. Hence, by applying the spectral
algorithm, the indices identifying the actors together with the teacher should correspond to the
children of a Q node. The other students, on the other hand, may be considered as isolated nodes
and so their indices, according to how the spectral seriation algorithm works, should identify the
children leaves of a P node.
As mentioned above, the spectral algorithm is applied to the similarity matrix associated with the matrix
A which describes the data and it constructs the PQ tree, displayed in Figure~\ref{pqtree} which encodes all the
admissible permutations of the nodes, i.e. the participants.

One of the admissible permutations, whose number is constrained by the particular structure of the tree,is
$$
\footnotesize{\texttt{p = [4 3 1 2 27 5 6 7 8 9 10 11 12 13 14 15 16 17 18 19 20 21 22 23 24 25 26]}}$$
and it represents the displayed frontier of the tree represented in Figure 1 that is one of all the possible
frontiers.
By observing the structure of the tree, it can be noted that the leaf nodes labelled as the students who
acted and the teacher, are isolated from the others which represent the observers.
Specifically, the nodes which represent the participants are children of a Q node that is characterized
by the fact that it represents both the right and the left orderings of its leaves labelled as the node in the
graph. The P node, child of the Q node, represents the equivalence class of all the permutations of its
leaves. Hence, all the constrained permutations of the nodes labelled as the participants in this phase,
i.e., the students who acted and the teacher, are the following

\begin{center}
\begin{tabular}{ccccc}
\hline\noalign{\smallskip}
4 & 3 & 1 & 2 & 27 \\
27 & 2 & 1 & 3 & 4 \\
4 & 3 & 2 & 1 & 27\\
27 & 1 & 2 & 3 & 4 \\
\noalign{\smallskip}\hline
\end{tabular}
\end{center}
\vspace{.5cm}

As aforementioned, the seriation problem consists in finding the best enumeration order of a set of units
so that elements with higher similarity are closer to each other in the resulting sequence. Under thisperspective, the obtained permutations allow us to make some considerations regarding the similarity
of the assigned roles based on the messages written in the chat. In fact, it can be observed that, in all
the admissible obtained permutations for the nodes that represent the active players in this episode,
nodes 1 and 2 (the Pest and the Blogger) are always close to each other and so are nodes 3 and 4 that
represent the Boss and the Promoter respectively. The teacher, labelled by node 27, is closely related to
the Pest and the Blogger as in all the resulting sequences these three nodes are placed together.
From the point of view of a didactic interpretation, the structure of the graph allows an analysis of the
behaviors of the student who actively participated with respect to the role assigned to each one of them
and with respect to the communication with the teacher. 
The expected behavior of the Promoter, whose
task is to animate the interaction between peers and any direct interaction with the uncle, is confirmed.
The student who assumes the role of the Pest, aimed at asking everyone questions, interacts mainly with
the Blogger on the mathematical question to be solved. We recall that this first episode is related to the
exploration of the problem from a mathematical point of view but also to the trial stage of the roles'
tasks. This may explain the particular structure of the graph related to the players which highlights a
similarity in the behavior of the roles of the Pest and the Blogger with respect to the other characters.
In this first phase, the teacher deliberately plays an interlocutory role with respect to both the
management of the group and the resolution of the mathematical task.
In the following, we analyse the data coming from the chat of the 5 groups of students who play as observers in the same analysed episode.
Each of these groups observed and communicate within itself simultaneously with the group of the actors who in the first episode had to ``explor'' the problem enclosed in the story. Note that students of each group interacted among themselfs and in this phase of the activity did not interact with the other observers nor with the teacher.
We considered the chat written simultaneously to the actors' one in order to study the behaviour of the students who observe the ones who play by comparing the PQ-trees of each observers' group.

The matrix coming from the analisys of the chat is constructed as mentioned above and has the following form:
$$
A=\left(\begin{array}{cccccccccccccccccccccccc} 1 & 0 & 0 & 0 & 0 & 0 & 0 & 0 & 0 & 0 & 0 & 0 & 0 & 0 & 0 & 0 & 0 & 0 & 0 & 0 & 16 & 0 & 2 & 0\\ 0 & 1 & 0 & 0 & 0 & 0 & 0 & 0 & 0 & 0 & 0 & 0 & 0 & 0 & 0 & 0 & 0 & 0 & 0 & 0 & 5 & 1 & 0 & 0\\ 0 & 0 & 1 & 0 & 0 & 0 & 0 & 0 & 0 & 0 & 0 & 0 & 0 & 0 & 0 & 0 & 0 & 0 & 0 & 0 & 14 & 3 & 0 & 1\\ 0 & 0 & 0 & 1 & 0 & 0 & 0 & 0 & 0 & 0 & 0 & 0 & 0 & 0 & 0 & 0 & 0 & 0 & 0 & 0 & 7 & 2 & 2 & 1\\ 0 & 0 & 0 & 0 & 1 & 0 & 0 & 0 & 0 & 0 & 0 & 0 & 0 & 0 & 0 & 0 & 0 & 0 & 0 & 0 & 10 & 0 & 0 & 0\\ 0 & 0 & 0 & 0 & 0 & 1 & 0 & 0 & 0 & 0 & 0 & 0 & 0 & 0 & 0 & 0 & 0 & 0 & 0 & 0 & 9 & 1 & 2 & 0\\ 0 & 0 & 0 & 0 & 0 & 0 & 1 & 0 & 0 & 0 & 0 & 0 & 0 & 0 & 0 & 0 & 0 & 0 & 0 & 0 & 9 & 0 & 2 & 0\\ 0 & 0 & 0 & 0 & 0 & 0 & 0 & 1 & 0 & 0 & 0 & 0 & 0 & 0 & 0 & 0 & 0 & 0 & 0 & 0 & 8 & 0 & 0 & 0\\ 0 & 0 & 0 & 0 & 0 & 0 & 0 & 0 & 1 & 0 & 0 & 0 & 0 & 0 & 0 & 0 & 0 & 0 & 0 & 0 & 1 & 0 & 0 & 0\\ 0 & 0 & 0 & 0 & 0 & 0 & 0 & 0 & 0 & 1 & 0 & 0 & 0 & 0 & 0 & 0 & 0 & 0 & 0 & 0 & 12 & 2 & 1 & 0\\ 0 & 0 & 0 & 0 & 0 & 0 & 0 & 0 & 0 & 0 & 1 & 0 & 0 & 0 & 0 & 0 & 0 & 0 & 0 & 0 & 9 & 1 & 0 & 0\\ 0 & 0 & 0 & 0 & 0 & 0 & 0 & 0 & 0 & 0 & 0 & 1 & 0 & 0 & 0 & 0 & 0 & 0 & 0 & 0 & 9 & 0 & 1 & 0\\ 0 & 0 & 0 & 0 & 0 & 0 & 0 & 0 & 0 & 0 & 0 & 0 & 1 & 0 & 0 & 0 & 0 & 0 & 0 & 0 & 3 & 10 & 8 & 8\\ 0 & 0 & 0 & 0 & 0 & 0 & 0 & 0 & 0 & 0 & 0 & 0 & 0 & 0 & 0 & 0 & 0 & 0 & 0 & 0 & 0 & 0 & 0 & 0\\ 0 & 0 & 0 & 0 & 0 & 0 & 0 & 0 & 0 & 0 & 0 & 0 & 0 & 1 & 0 & 0 & 0 & 0 & 0 & 0 & 2 & 12 & 5 & 3\\ 0 & 0 & 0 & 0 & 0 & 0 & 0 & 0 & 0 & 0 & 0 & 0 & 0 & 0 & 1 & 0 & 0 & 0 & 0 & 0 & 7 & 17 & 10 & 3\\ 0 & 0 & 0 & 0 & 0 & 0 & 0 & 0 & 0 & 0 & 0 & 0 & 0 & 0 & 0 & 0 & 0 & 0 & 0 & 0 & 0 & 0 & 0 & 0\\ 0 & 0 & 0 & 0 & 0 & 0 & 0 & 0 & 0 & 0 & 0 & 0 & 0 & 0 & 0 & 1 & 0 & 0 & 0 & 0 & 3 & 4 & 1 & 4\\ 0 & 0 & 0 & 0 & 0 & 0 & 0 & 0 & 0 & 0 & 0 & 0 & 0 & 0 & 0 & 0 & 1 & 0 & 0 & 0 & 1 & 12 & 9 & 3\\ 0 & 0 & 0 & 0 & 0 & 0 & 0 & 0 & 0 & 0 & 0 & 0 & 0 & 0 & 0 & 0 & 1 & 0 & 0 & 0 & 0 & 0 & 0 & 0\\ 0 & 0 & 0 & 0 & 0 & 0 & 0 & 0 & 0 & 0 & 0 & 0 & 0 & 0 & 0 & 0 & 0 & 1 & 0 & 0 & 1 & 8 & 9 & 4\\ 0 & 0 & 0 & 0 & 0 & 0 & 0 & 0 & 0 & 0 & 0 & 0 & 0 & 0 & 0 & 0 & 0 & 0 & 0 & 0 & 0 & 0 & 0 & 0\\ 0 & 0 & 0 & 0 & 0 & 0 & 0 & 0 & 0 & 0 & 0 & 0 & 0 & 0 & 0 & 0 & 0 & 0 & 1 & 0 & 2 & 10 & 15 & 2\\ 0 & 0 & 0 & 0 & 0 & 0 & 0 & 0 & 0 & 0 & 0 & 0 & 0 & 0 & 0 & 0 & 0 & 0 & 0 & 1 & 1 & 14 & 11 & 5\\ 0 & 0 & 0 & 0 & 0 & 0 & 0 & 0 & 0 & 0 & 0 & 0 & 0 & 0 & 0 & 0 & 0 & 0 & 0 & 0 & 0 & 0 & 0 & 0 \end{array}\right);
$$

In particular, the rows represent the students who are observers. All the 26 students are divided into groups $g_i$; the members of $g_1$ are the actors while the observers are in the other groups. Groups $g_i$  for $i = 2,3,4$ consist of four students, in group $g_5$ there were six students therefore 2 students play the PESTE role and two perform the BOSS. The last group was composed of six students, two of them perform the PESTE and two had the PROMOTER's role.
Given the composition of the observers' groups, the rows of the matrix are 25. 
There are 24 columns; for each of the 5 groups of observers there are 4 columns relating to the roles of Pest,
Blogger, Boss and Promoter. Therefore the first 20 columns represent the roles and, as in the case
of the actors' chat, the remaining columns represent the different types of interactions
minutes. In particular, column 21 contains the main interventions, i.e., those relating to the observed role
corresponding to that assumed by the considered observer pupil. In column 22 the
secondary / social interventions are contained; elements $a_{i,23}$ are
general interventions on the behavior of the group of actors or on a role
different from that assumed by the $i$-th observer. Finally, the last column contains
interventions relating to the mathematical problem or to solution proposed by the students in the chat 
group that observers discuss.

We report and analyze the peer interactions of the observers in each group, in order to construct the PQtree coming from processing the similarity matrix costructed from the binary data matrix. 
Note that the latter does not take into account the frequencies but is obtained by replacing each non-zero entry of $A$ by 1. 

For the 5 observers' groups the output of the algorithm is given by a tree with only a Q node, two trees rooted by a P node and the trees of groups 4 and 6 have a pattern similar to the actors' one with a mized P node and Q node structure.

Riportiamo alcuni estratti e analizziamo le interazioni delle chat di ciascun gruppo, al fine di identificare se alla analogia di struttura dell'albero corrisponda una analogia di comportamento tra gruppi di osservatori o con il gruppo degli attori.

In the following we report some excerpts and then analyze the interactions of each group, in order to identify whether the analogy of the structure of the tree corresponds to an analogy of behavior between groups of observers or with the group of actors.

By analyzing the observers' similarity matrix the obtained permutations is given by the sequence of nodes

$$
\footnotesize{\texttt{p = [5 8 9 1 7 12 2 11 6 10 3 4 13 15 16 18 21 23 24 19 20 14 17 22 25]}}.$$

In the following, $B_k$ indicates the data matrix relating to the group of observers $ g_k $ with $ k = 2, \ldots, 6 $. 

\medskip
\textbf{Group $g_2$.}

\begin{figure}
\begin{center}
\includegraphics[width=.5\textwidth]{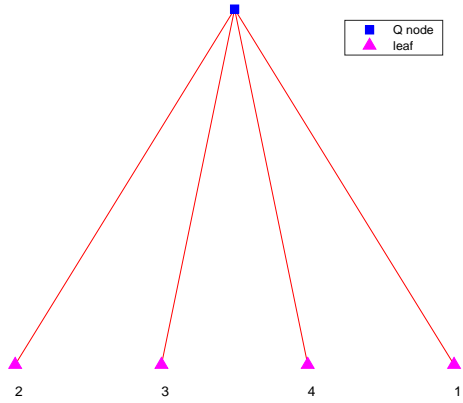}
\caption{PQ-tree corresponding to the data of group  $g_2$.}
\label{g2}
\end{center}
\end{figure}

The data matrix related to observers' group $g_2$ and the similarity matrix of the corresponding ${0,1}$ data matrix are respectively

$$
B_2 =\left(\begin{array}{cccccccc} 1 & 0 & 0 & 0 & 16 & 0 & 2 & 0\\ 0 & 1 & 0 & 0 & 5 & 1 & 0 & 0\\ 0 & 0 & 1 & 0 & 14 & 3 & 0 & 1\\ 0 & 0 & 0 & 1 & 7 & 2 & 2 & 1 \end{array}\right), \quad     S_2\left(\begin{array}{cccc} 3 & 1 & 1 & 2\\ 1 & 3 & 2 & 2\\ 1 & 2 & 4 & 3\\ 2 & 2 & 3 & 5 \end{array}\right)
$$

and one of the admissible permutations encoded in the PQtree obtained through the seriation algorithm is
\begin{verbatim}
 p2 = [2 3 4 1]
\end{verbatim}
that is the frontier of the PQtree, depicted in Figure~\ref{g2}, with a Q node as a root. 

We observed that in this group only the Boss and the Promoter made considerations on the problem, like the following
\medskip

 Promoter 12:50: \textit{according to the Boss, the Promoter did not write the correct formula};

 Boss 12:51: \textit{the Boss suggests another formula}

\medskip

The contingency analysis shows that in this group there were few social interactions and each member commented on the same interpreted role. In particular the Pest posted many comments related to the actor Pest.
Group $g_2$ is the only group that has a compact and uniform behavior with respect to the types
of interactions we are analyzing. All the other observers' groups are characterized by behaviors more or less distinct within the roles or type of interventions.

\medskip
\textbf{Group $g_4$.}

\begin{figure}
\begin{center}
\includegraphics[width=.5\textwidth]{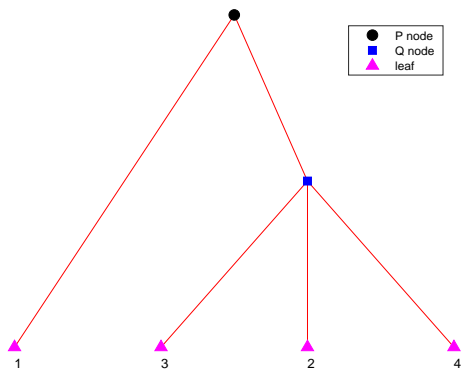}
\caption{PQ-tree corrispondente alla matrice dei dati $B_4$.}
\label{g4}
\end{center}
\end{figure}

The matrix constructed from the analysis of the data in the chat of group $g_4$ and the similarity matrix constructed from its binary version are respectively

$$
B_4=\left(\begin{array}{cccccccc} 1 & 0 & 0 & 0 & 1 & 0 & 0 & 0\\ 0 & 1 & 0 & 0 & 12 & 2 & 1 & 0\\ 0 & 0 & 1 & 0 & 9 & 1 & 0 & 0\\ 0 & 0 & 0 & 1 & 9 & 0 & 1 & 0 \end{array}\right), \quad     S_4\left(\begin{array}{cccc} 2 & 1 & 1 & 1\\ 1 & 4 & 2 & 2\\ 1 & 2 & 3 & 1\\ 1 & 2 & 1 & 3 \end{array}\right)
$$

and one of the admissible permutations encoded in the obtained tree is the sequence
\begin{verbatim}
p4 = [1 3 2 4]. 
\end{verbatim}
that is the frontier of the PQtree in Figure~\ref{g4}.
This PQtree has a P node as root and has two branching that identify different behaviours
between node 1 (the Pest) and the other members. In this group the Pest has a marginal involvement at the end of this episode concerning a consideration of how the Pest-actor is playing

\medskip
Pest 13:10: \textit{the Pest is performing very well the role.}

%

\medskip
\textbf{Group $g_6$.}

The data matrix coming form the chat of group $g_6$ and the related similarity matrix of the binary version are respectively
\medskip

$$
B_6 =\left(\begin{array}{cccccccc} 1 & 0 & 0 & 0 & 1 & 11 & 8 & 3\\ 1 & 0 & 0 & 0 & 0 & 1 & 1 & 0\\ 0 & 1 & 0 & 0 & 1 & 8 & 9 & 4\\ 0 & 0 & 1 & 0 & 2 & 10 & 15 & 2\\ 0 & 0 & 0 & 1 & 1 & 14 & 11 & 5 \end{array}\right), \quad S_6 =     \left(\begin{array}{ccccc} 5 & 3 & 4 & 4 & 4\\ 3 & 3 & 2 & 2 & 2\\ 4 & 2 & 5 & 4 & 4\\ 4 & 2 & 4 & 5 & 4\\ 4 & 2 & 4 & 4 & 5 \end{array}\right).
$$

\medskip

\begin{figure}
\begin{center}
\includegraphics[width=.4\textwidth]{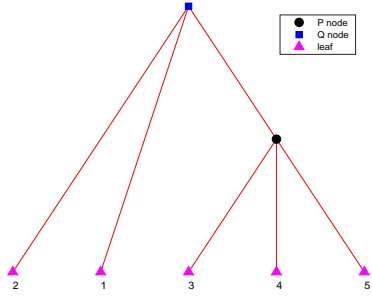}
\caption{PQ-tree related to the data matrix$B_6$.}
\label{g6}
\end{center}
\end{figure}

The PQ-tree obtained after applying the seriation algorithm to the simularity matrix constructed from the binary version of $B_6$ is depicted  in Figure~\ref{g6}.


One of the admissible permutations, that is the frontier of the PQtree in Figure~\ref{g6}, is the sequence
\begin{verbatim}
p6 = [2 1 3 4 5].
\end{verbatim}
\medskip
\textbf{Group $g_5$.}

\begin{figure}
\begin{center}
\includegraphics[width=.5\textwidth]{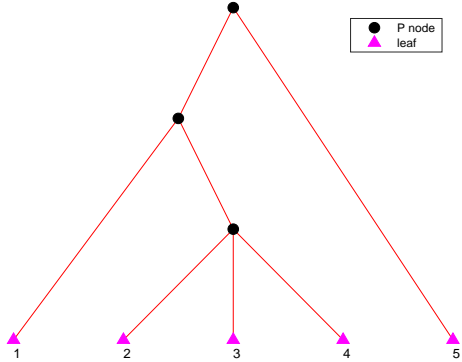}
\caption{PQ-tree related to the analysis of the matrix $B_5$.}
\label{g5}
\end{center}
\end{figure}

The data matrix related to the chat of group $g_5$ and the similarity matrix are 
$$
B_5 = \left(\begin{array}{cccccccc} 1 & 0 & 0 & 0 & 1 & 0 & 0 & 0\\ 0 & 1 & 0 & 0 & 2 & 12 & 5 & 3\\ 0 & 0 & 1 & 0 & 7 & 17 & 10 & 3\\ 0 & 0 & 0 & 1 & 3 & 4 & 1 & 4\\ 0 & 0 & 0 & 0 & 0 & 0 & 0 & 0 \end{array}\right), \quad S_5=     \left(\begin{array}{ccccc} 2 & 1 & 1 & 1 & 0\\ 1 & 5 & 4 & 4 & 0\\ 1 & 4 & 5 & 4 & 0\\ 1 & 4 & 4 & 5 & 0\\ 0 & 0 & 0 & 0 & 0 \end{array}\right)
$$
and one of the admissible permutations, that is the frontier of the tree depicted in Figure~\ref{g5} is the sequence 
\begin{verbatim}
p5 = [1 2 3 4 5]. 
\end{verbatim}

The PQtree in Figure~\ref{g5} represents in a compact way the following 24 permutations
\begin{center}
\begin{tabular}{ccccc}
\hline\noalign{\smallskip}
5 & 4 & 3 & 2 & 1\\
5  &   4  &   2  &   3  &   1\\
5  &   3  &   4  &   2  &   1\\
5  &   3  &   2  &   4  &   1\\
5  &   2  &   4  &   3  &   1\\
5  &  2   &  3  &   4  &   1\\
5  &   1  &   4  &   3   &  2\\
5  &   1  &   4 &    2  &   3\\
5  &   1   &  3 &    4 &    2\\
5  &   1   &  3  &   2  &   4\\
5  &   1  &   2  &   4  &   3\\
5  &   1  &   2  &   3  &   4\\
4  &   3   &  2  &   1  &   5\\
4  &   2   &  3   &  1  &   5\\
3  &   4   &  2   &  1  &   5\\
3  &   2   &  4  &   1  &   5\\
2  &   4   &  3  &   1  &   5\\
2   &  3 &    4  &   1 &    5\\
1   &  4  &   3 &    2  &   5\\
1  &   4   &  2  &   3  &   5\\
1    & 3  &   4  &   2  &   5\\
1  &   3   &  2   &  4   &  5\\
1   &  2  &   4  &   3  &   5\\
1   &  2  &   3   &  4 &    5\\
\noalign{\smallskip}\hline
\end{tabular}
\end{center}

\medskip
\textbf{Group $g_3$.}

\begin{figure}
\begin{center}
\includegraphics[width=.5\textwidth]{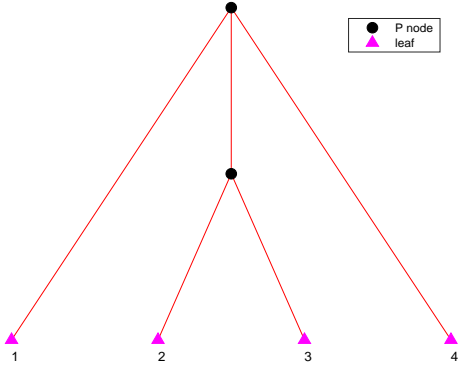}
\caption{PQ-tree related to the analysis of group $g_3$.}
\label{g3}
\end{center}
\end{figure}

The third observers' group is characterized by similar behavior of all members who did not
make comments on the problem. 
The data matrix coming form the chat of the observers in group $g_3$ and the similarity matrix constructed form the binary data matrix are respectively
$$
B_3 = \left(\begin{array}{cccccccc} 1 & 0 & 0 & 0 & 10 & 0 & 0 & 0\\ 0 & 1 & 0 & 0 & 9 & 1 & 2 & 0\\ 0 & 0 & 1 & 0 & 9 & 0 & 2 & 0\\ 0 & 0 & 0 & 1 & 8 & 0 & 0 & 0 \end{array}\right), \quad S_3 = \left(\begin{array}{cccc} 2 & 1 & 1 & 1\\ 1 & 4 & 2 & 1\\ 1 & 2 & 3 & 1\\ 1 & 1 & 1 & 2 \end{array}\right)
$$
one of the possible permutation is the frontier of the PQtree depicted in Figure~\ref{g3}, that is the sequence
\begin{verbatim}
p3 = [1 2 3 4]. 
\end{verbatim}

The particular structure of the tree encode the following 12 admissible permutations

\begin{center}
\begin{tabular}{cccc}
\hline\noalign{\smallskip}
4 & 3 & 1 & 2\\
4 & 2 & 3 & 1\\
4 & 1 & 3 & 2\\
4 & 1 & 2 & 3\\
3 & 2 & 4 & 1\\
2 & 3 & 4 & 1\\
3 & 2 & 1 & 4\\
2 & 3 & 1 & 4\\
1 & 4 & 3 & 2\\
1 & 4 & 2 & 3\\
1 & 3 & 2 & 4\\
1 & 2 & 3 & 4\\
\noalign{\smallskip}\hline
\end{tabular}
\end{center}
\vspace{.5cm}

Reading the chat we observed that the Pest (node 1) and the Promoter (node 4) made 
considerations only on the role of their correspondent-actors and in fact in the tree in Figure~\ref{g3} are separated, as it should have been, from the Boss and the Blogger who instead intervened on the actors who play other roles. This can be observed in the following posts reported here as examples:

Promoter $g_3$ 12:30: \textit{the Promoter takes the floor allowing to...}

Promoter $g_3$ 12:31: \textit{the Promoter gives proof of what he/she is saying}

Peste $g_3$ 12:31: \textit{the Pest is participating a lot asking many questions}

Blogger $g_3$ 12:32: \textit{the Blogger is participating a lot...}

Boss $g_3$ 12:33: \textit{the Boss is coordinating the group very well...}

Reading the observers' interventions, it was noted that groups $g_5$ and $g_6$ are more oriented to comment in general on the behavior of the actors and not on the specific role that is the same one played by them. Moreover, rather than analyzing the behavior of the actors who fulfill their corresponding role, the observers of these two groups interact by making interventions also on the mathematical problem.


\section{Conclusions}
In this manuscript, we presented the preliminary results of an experimental work, on digital
storytelling in Mathematics, aimed at promoting critical mathematical thinking among students.
The analyzed activity is the first episode of a story that develops in 5 episodes and involves the
participation of students from a second-grade class who take on both the role of active protagonists
and that of observers of the actors during the problem solving activity. In particular, the
communication between the participants who interacted on the chat of the Moodle e-learning
platform was examined from a didactic and numerical modeling point of view.
The case study analyzed above, allows us to validate the spectral algorithm and the accompanying
Matlab toolbox originally implemented for solving the seriation problem in archaeology where the
latter appeared first. In particular, the outcomes obtained for the active members of the study and
the comments at the end of Section 3, show that the spectral algorithm gives a reasonable result in
terms of admissible permutations for the indices of the nodes depicting the assigned roles.
Although the size of the examined sample is not very large, the structure of the algorithm allows
us to confirm that also for wider dimensions, the group of the active participants should be
separated from the observers as it identifies a connected component in the underneath graph.
Moreover, from a numerical point of view, the size of the similarity matrix does not constitute a
problem since the spectral algorithm is based on the computation of an eigenvector and for thispurpose the toolbox includes the possibility to choose between a small-scale and a large-scale
implementation and, for improving the performance in case of a large-scale problem, a parallel
version of the method is also provided.
In this work, we observed how the complex networks analysis in a DIST-M setting may be worthy both
as a tool for validating the mentioned algorithm and as a potential instrument for allowing the
identification of the behaviors, whether didactic or interpersonal, one wants to examine.
Further development of this research may include the analogous investigation of the episodes
subsequent to the first analyzed in this manuscript, in order to put emphasis on potential dissimilarities
in the behavior of the students when they assume a different role.
%

\end{document}